# IDENTIFIABILITY OF PARAMETERS IN LATENT STRUCTURE MODELS WITH MANY OBSERVED VARIABLES[1]

By Elizabeth S. Allman[2], Catherine Matias and John A. Rhodes[2]

*University of Alaska Fairbanks, Université d'Évry val d'Essonne, CNRS and University of Alaska, Fairbanks*

While hidden class models of various types arise in many statistical applications, it is often difficult to establish the identifiability of their parameters. Focusing on models in which there is some structure of independence of some of the observed variables conditioned on hidden ones, we demonstrate a general approach for establishing identifiability utilizing algebraic arguments. A theorem of J. Kruskal for a simple latent-class model with finite state space lies at the core of our results, though we apply it to a diverse set of models. These include mixtures of both finite and nonparametric product distributions, hidden Markov models and random graph mixture models, and lead to a number of new results and improvements to old ones.

In the parametric setting, this approach indicates that for such models, the classical definition of identifiability is typically too strong. Instead generic identifiability holds, which implies that the set of nonidentifiable parameters has measure zero, so that parameter inference is still meaningful. In particular, this sheds light on the properties of finite mixtures of Bernoulli products, which have been used for decades despite being known to have nonidentifiable parameters. In the nonparametric setting, we again obtain identifiability only when certain restrictions are placed on the distributions that are mixed, but we explicitly describe the conditions.

**1. Introduction.** Statistical models incorporating latent variables are widely used to model heterogeneity within datasets, via a hidden structure. How-

Received September 2008; revised January 2009.

[1]Supported by the Isaac Newton Institute and the Statistical and Applied Mathematical Sciences Institute.

[2]Supported by NSF Grant DMS-07-14830.

ever, the fundamental theoretical question of the identifiability of parameters of such models can be difficult to address. For specific models it is even known that certain parameter values lead to nonidentifiability, while empirically, the model appears to be well behaved for most values. Thus parameter inference procedures may still be performed, even though theoretical justification of their consistency is still lacking. In some cases (e.g., hidden Markov models [39]), it has been formally established that generic choices of parameters are identifiable, which means that only a subset of parameters of measure zero may not be identifiable.

In this work, we consider a number of such variable models, all of which exhibit a conditional independence structure, in which (some of) the observed variables are independent when conditioned on the unobserved ones. In particular, we investigate:

1. finite mixtures of products of finite measures, where the mixing parameters are unknown (including finite mixtures of multivariate Bernoulli distributions), also called *latent-class models* in the literature;
2. finite mixtures of products of nonparametric measures, again with unknown mixing parameters;
3. discrete hidden Markov models;
4. a random graph mixture model, in which the probability of the presence of an edge is determined by the hidden states of the vertices it joins.

We show how a fundamental algebraic result of Kruskal [29, 30] on 3-way tables can be used to derive identifiability results for all of these models. While Kruskal's work is focused on only 3 variates, each with finite state spaces, we use it to obtain new identifiability results for mixtures with more variates (point 1, above), whether discrete or continuous (point 2). For hidden Markov models (point 3), with their more elaborate dependency structure, Kruskal's work allows us to easily recover some known results on identifiability that were originally approached with other tools, and to strengthen them in certain aspects. For the random graph mixture model (point 4), in which the presence/absence of each edge is independent conditioned on the states of all vertices, we obtain new identifiability results via this method, again by focusing on the model's essential conditional independence structure.

While we establish the validity of many identifiability statements not previously known, the major contribution of this paper lies as much in the method of analysis it introduces. By relating a diverse collection of models to Kruskal's work, we indicate the applicability of this method of establishing identifiability to a variety of models with appropriate conditional independence structure. Although our example of applying Kruskal's work to a complicated model such as the random graph requires substantial additional algebraic arguments tied to the details of the model, it illustrates well that the essential insight can be a valuable one.

Finally, we note that in establishing identifiability of the parameters of a model, this method clearly indicates one must allow for the possibility of certain "exceptional" choices of parameter values which are not identifiable. However, as these exceptional values can be characterized through algebraic conditions, one may deduce that they are of measure zero within the parameter space (in the finite-dimensional case). Since "generic" parameters are identifiable, one is unlikely to face identifiability problems in performing inference. Thus generic identifiability of the parameters of a model is generally sufficient for data analysis purposes. Although the notion of identifiability of parameters off a set of measure zero is not a new one, neither the usefulness of this notion nor its algebraic origins seem to have been widely recognized.

**2. Background.** Latent structure models form a very large class of models including, for instance, finite univariate or multivariate mixtures [34], hidden Markov models [5, 16] and nonparametric mixtures [33].

General formulations of the identification problem were made by several authors, and pioneering works may be found in [27, 28]. The study of identifiability proceeds from a hypothetical exact knowledge of the distribution of observed variables and asks whether one may, in principle, recover the parameters. Thus identification problems are not problems of statistical inference in a strict sense. However, since nonidentifiable parameters cannot be consistently estimated, identifiability is a prerequisite of statistical parameter inference.

In the following, we are interested in models defined by a family $\mathcal{M}(\Theta) = \{\mathbb{P}_\theta, \theta \in \Theta\}$ of probability distributions on some space $\Omega$, with parameter space $\Theta$ (not necessarily finite dimensional). The classical definition of identifiability, which we will refer to as *strict identifiability*, requires that for any two different values $\theta \neq \theta'$ in $\Theta$, the corresponding probability distributions $\mathbb{P}_\theta$ and $\mathbb{P}_{\theta'}$ are different. This is equivalent to injectivity of the model's parameterization map $\Psi$, which takes values in $\mathcal{M}_1(\Omega)$, the set of probability measures on $\Omega$, and is defined by $\Psi(\theta) = \mathbb{P}_\theta$.

In many cases, the above map will not be strictly injective. For instance, it is well known that in models with discrete hidden variables (such as finite mixtures or discrete hidden Markov models), the latent classes can be freely relabeled without changing the distribution of the observations, a phenomenon known as "label swapping." In this sense, the above map is always at least $r!$-to-one, where $r$ is the number of classes in the model. However, this does not prevent the statistician from inferring the parameters of these models. Indeed, parameter identifiability up to a permutation on the class labels (which we henceforth consider as a type of strict identifiability), is largely enough for practical use, at least in a maximum likelihood setting. Note that the label swapping issue may cause major problems in a Bayesian framework (see, for instance, [34], Section 4.9).

A related concept of *local identifiability* only requires the parameter to be unique in small neighborhoods in the parameter space. For parametric models (i.e., when the parameter space is finite dimensional), with some regularity conditions, there is an equivalence between local identifiability of the parameters and nonsingularity of the information matrix [40]. When an iterative procedure is used to approximate an estimator of the parameter, different initializations can help to detect multiple solutions of the estimation problem. This often corresponds to the existence of multiple parameter values giving rise to the same distribution. However, the validity of such procedures relies on knowing that the parameterization map is, at most, finite-to-one, and a precise characterization of the value of $k$ such that it is a $k$-to-one map would be most useful.

Thus knowledge that the parameterization map is finite-to-one might be too weak a result from a statistical perspective on identifiability. Moreover, we argue in the following that infinite-to-one maps might not be problematic, as long as they are *generic* $k$-to-one maps for known finite $k$.

While all our results are proved relying on the same underlying tool, they must be expressed differently in the parametric framework (including the finite case) and in the nonparametric one.

*The parametric framework.* While the focus on one-to-one or $k$-to-one parameterization maps is well suited for most of the classical models encountered in the literature, it is inadequate in some important cases. For instance, it is well known that finite mixtures of Bernoulli products are not identifiable [23], even up to a relabeling of latent classes. However, these distributions are widely used to model data when many binary variables are observed from individuals belonging to different unknown populations, and parameter estimation procedures are performed in this context. For instance, these models may be used in numerical identification of bacteria (see [23] and the references therein). Statisticians are aware of this apparent contradiction; the title of the article [6], *practical identifiability of finite mixtures of multivariate Bernoulli distributions*, indicates the need to reconcile nonidentifiability and validity of inference procedures, and clearly indicates that the strict notion of identifiability is not useful in this specific context. We establish that parameters of finite mixtures of multivariate Bernoulli distributions (with a fixed number of components) are in fact generically identifiable (see Section 5).

Here, "generic" is used in the sense of algebraic geometry, as will be defined in the subsection on algebraic terminology below. Most importantly, it implies that the set of points for which identifiability does not hold has measure zero. In this sense, any observed data set has probability one of being drawn from a distribution with identifiable parameters.

Understanding when generic identifiability holds, even in the case of finite measures, can be mathematically difficult. There are well-known examples of latent-class models in which the parameterization map is in fact infinite-to-one, for reasons that are not immediately obvious. For instance, Goodman [22] describes a 3-class model with four manifest binary variables and thus a parameter space of dimension $3(4) + 2 = 14$. Though the distributions resulting from this model lie in a space of dimension $2^4 - 1 = 15$, the image of the parameterization map has dimension only 13. From a statistical point of view, this results in nonidentifiability.

An important observation that underlies our investigations is that many finite space models (e.g., latent-class models, hidden Markov models) involve parameterization maps which are polynomial in the scalar parameters. Thus statistical models have recently been studied by algebraic geometers [19, 37]. Even in the more general case of distributions belonging to an exponential family, which lead to analytic but nonpolynomial maps, it is possible to use perspectives from algebraic geometry (see, for instance, [2, 3, 13]). Algebraic geometers use terminology rather different from the statistical language, for instance, they describe the image of the parameterization map of a simple latent-class model as a higher secant variety of a Segre variety. When the dimension of this variety is less than expected, as in the example of Goodman above, the variety is termed defective, and one may conclude the parameterization map is generically infinite-to-one. Recent works such as [1, 7, 8] have made much progress in determining when defects occur.

However, as pointed out by Elmore, Hall and Neeman [14], focusing on dimension is not sufficient for a complete understanding of the identifiability question. Indeed, even if the dimensions of the parameter space and the image match, the parameterization might be a generically $k$-to-one map, and the finite number $k$ cannot be characterized by using dimension approaches. For example, consider latent-class models, assuming the number $r$ of classes is known. In this context, even though the dimensions agree, we might have a generically $k$-to-one map with $k > r!$. (Recall that $r!$ corresponds to the number of points which are equivalent by permutating label classes.)

This possibility was already raised in the context of psychological studies by Kruskal [29], whose work in [30] provides a strong result ensuring generic $r!$-to-oneness of the parameterization map for latent $r$-class models under certain conditions. Kruskal's work, however, is focused on models with only 3 observed variables, or, in other terms, on secant varieties of Segre products with 3 factors, or on 3-way arrays. While the connection of Kruskal's work to the algebraic geometry literature seems to have been overlooked, the nature of his result is highly algebraic.

Although [14] is ultimately directed at understanding nonparametric mixtures, Elmore, Hall, and Neeman address the question of identifiability of the parameterization for latent-class models with many binary observed variables (i.e., for secant varieties of Segre products of projective lines with many factors, or on $2 \times 2 \times \cdots \times 2$ tables). These are just the mixtures of Bernoulli products referred to above, though the authors never introduce that terminology. Using algebraic methods, they show that with sufficiently many observed variables, the image of the parameterization map is birationally equivalent to a symmetrization of the parameter space under the symmetric group $\Sigma_r$. Thus for sufficiently many observed variables, the parameterization map is generically $r!$-to-one. (Although the generic nature of the result is not made explicit, that is, however, all that one can deduce from a birational equivalence.) Their proof is constructive enough to give a numerical understanding of how many observed variables are sufficient, though this number's growth in $r$ is much larger than is necessary (see Corollary 5 and Theorem 8 for more details).

*The nonparametric framework.*  Nonparametric mixture models have received much attention recently [4, 25, 26]. They provide an interesting framework for modelling very general heterogeneous data. However, identifiability is a difficult and crucial issue in such a high-dimensional setting.

Using algebraic methods to study statistical models is most straightforward when state spaces are finite. One way of handling continuous random variables via an algebraic approach is to discretize the problem by binning the random variable into a finite number of sets. For instance, [11, 15, 26] developed cut points methods to transform multivariate continuous observations into binomial or multinomial random variables.

As already mentioned, Elmore, Hall and Neeman [14] consider a finite mixture of products of continuous distributions. By binning each continuous random variable $X$ to create a binary one, defined by the indicator function $1\{X \le t\}$, for some choice of $t$, they pass to a related finite model. But identification of a distribution is equivalent to identification of its cumulative distribution function (c.d.f.) $F(t) = \mathbb{P}(X \le t)$. Having addressed the question of identifiability of the parameters of a mixture of products of binary variables, they can thus argue for the identifiability of the parameters of the original continuous model, as they continue to do in [24]. However, because the authors are not explicit about the generic aspect of their results in [14], there are significant gaps in the formal justification of their claims. Moreover, the bounds they claim on the number of observed variables which ensure generic identifiability leave much room for improvement, as they point out.

*The general approach.*  Our theme in this work is the applicability of the fundamental result of Kruskal on 3-way arrays to a spectrum of models with latent structure. Though our approach is highly algebraic, it has little in common with that of [14, 24], for establishing that with sufficiently many observed variables, the parameterization map of $r$-latent-class models is either generically $r!$-to-one in the parametric case, or that it is exactly $r!$-to-one (under some conditions) in the nonparametric case. Our results apply not only to binary variables, but as easily to ones with more states, or even to continuous ones. In the case of binary variables (multivariate Bernoulli mixtures), we obtain a much lower upper bound for a sufficient number of variables to ensure generic identifiability (up to label swapping) than the one that can be deduced from [14], and, in fact, our bound gives the correct order of growth, $\log_2 r$. (The constant factor we obtain is, however, still unlikely to be optimal.)

While our first results are on the identifiability of finite mixtures (with a fixed number of components) of finite measure products, our method has further consequences for more sophisticated models with a latent structure. Our approach for such models with finite state spaces can be summarized very simply: we group the observed variables into 3 collections, and view the composite states of each collection as the states of a single clumped variable. We choose our collections so that they will be conditionally independent, given the states of some of the hidden variables. Viewing these hidden variables as a single composite one, the model reduces to a special instance of the model Kruskal studied. Thus Kruskal's result on 3-way tables can be applied, after a little work, to show that Kruskal's condition is satisfied. This might be done either by showing that the clumping process results in a sufficiently generic model (ensuring Kruskal's condition is automatically satisfied for generic parameters), or that explicit restrictions on the parameters ensure this clumping process satisfies Kruskal's condition. In more geometric terms, we embed a complicated finite model into a simple latent-class model with 3 observed variables, taking care to verify that the embedding does not end up in the small set for which Kruskal's result tells us nothing.

To take up the continuous random variables case, we simply bin the real-valued random variables into a partition of $\mathbb{R}$ into $\kappa$ intervals and apply the previous method to the new discretized random variables. As a consequence, we are able to prove that finite mixtures of nonparametric independent variates, with at least 3 variates, have identifiable parameters under a mild and explicit regularity condition. This is in sharp contrast not only with [14, 24], but also with works such as [26], where the components of the mixture are assumed to be independent but also identically distributed and [25], which dealt only with $r = 2$ groups (see Section 7 for more details).

We note that Kruskal's result has already been successfully used in phylogeny, to prove identifiability of certain models of evolution of biological sequences along a tree [3]. However, application of Kruskal's result is limited to hidden class models, or to other models with some conditional independence structure, which have at least 3 observed variates. Kruskal's theorem can sometimes be used for models with many hidden variables, by considering a clumped latent variable $Z = (Z_1, \ldots, Z_n)$. We give two examples of such a use for models presenting a dependency structure on the observations, namely hidden Markov models (Section 6.1) and mixture models for random graphs (Section 6.2). For hidden Markov models, we recover many known results, and improve on some of them. For the random graph mixture model, we establish identifiability for the first time. Note that in all these applications we always assume the number of latent classes is known, which is crucial in using Kruskal's approach. Identification of the number of classes is an important issue that we do not consider here.

*Algebraic terminology.* Polynomials play an important role throughout our arguments, so we introduce some basic terminology and facts from algebraic geometry that we need. For a more thorough but accessible introduction to the field, we recommend [10].

An *algebraic variety* $V$ is defined as the simultaneous zero-set of a finite collection of multivariate polynomials $\{f_i\}_{i=1}^n \subset \mathbb{C}[x_1, x_2, \ldots, x_k]$,

$$V = V(f_1, \ldots, f_n) = \{\mathbf{a} \in \mathbb{C}^k | f_i(\mathbf{a}) = 0, 1 \leq i \leq n\}.$$

A variety is all of $\mathbb{C}^k$ only when all $f_i$ are 0; otherwise, a variety is called a *proper subvariety* and must be of dimension less than $k$, and, hence, of Lebesgue measure 0 in $\mathbb{C}^k$. Analogous statements hold if we replace $\mathbb{C}^k$ by $\mathbb{R}^k$, or even by any subset $\Theta \subseteq \mathbb{R}^k$ containing an open $k$-dimensional ball. This last possibility is of course most relevant for the statistical models of interest to us, since the parameter space is naturally identified with a full-dimensional subset of $[0, 1]^L$ for some $L$ (see Section 3 for more details). Intersections of algebraic varieties are algebraic varieties as they are the simultaneous zero-set of the unions of the original sets of polynomials. Finite unions of varieties are also varieties, since if sets $S_1$ and $S_2$ define varieties, then $\{fg | f \in S_1, g \in S_2\}$ defines their union.

Given a set $\Theta \subseteq \mathbb{R}^k$ of full dimension, we will often need to say some property holds for all points in $\Theta$, except possibly for those on some proper subvariety $\Theta \cap V(f_1, \ldots, f_n)$. We express this succinctly by saying the property holds *generically* on $\Theta$. We emphasize that the set of exceptional points of $\Theta$, where the property need not hold, is thus necessarily of Lebesgue measure zero.

In studying parametric models, $\Theta$ is typically taken to be the parameter space for the model, so that a claim of *generic identifiability* of model parameters means that all nonidentifiable parameter choices lie within a proper subvariety, and thus form a set of Lebesgue measure zero. While we do not always explicitly characterize the subvariety in statements of theorems, one could do so by careful consideration of our proofs.

In a nonparametric context, where algebraic terminology appropriate to the finite-dimensional setting is inappropriate, we avoid the use of the term "generic." Instead, we always give explicit characterizations of those parameter choices which may not be identifiable.

*Roadmap.* We first present finite mixtures of finite measure products with a conditional independence structure (or latent-class models) in Section 3. Then, Kruskal's result and consequences are presented in Section 4. Direct consequences on the identifiability of the parameters of finite mixtures of finite measure products appear in Section 5. More complicated dependent variables models, including hidden Markov models and a random graph mixture model, are studied in Section 6. In Section 7, we consider mixtures of nonparametric distributions, analogous to the finite ones considered earlier. All proofs are postponed to Section 8.

**3. Finite mixtures of finite measure products.** Consider a vector of observed random variables $\{X_j\}_{1 \leq j \leq p}$ where $X_j$ has finite state space with cardinality $\kappa_j$. Note that these variables are not assumed to be i.i.d. nor to have the same state space. To model the distribution of these variables, we use a latent (unobserved) random variable $Z$ with values in $\{1, \ldots, r\}$, where $r$ is assumed to be known. We interpret $Z$ as denoting an unobservable class, and assume that conditional on $Z$, the $X_j$'s are independent random variables. The probability distribution of $Z$ is given by the vector $\boldsymbol{\pi} = (\pi_i) \in (0,1)^r$ with $\sum \pi_i = 1$. Moreover, the probability distribution of $X_j$ conditional on $Z = i$ is specified by a vector $\mathbf{p}_{ij} \in [0,1]^{\kappa_j}$. We use the notation $\mathbf{p}_{ij}(l)$ for the $l$th coordinate of this vector $(1 \leq l \leq \kappa_j)$. Thus we have $\sum_l \mathbf{p}_{ij}(l) = 1$.

For each class $i$, the joint distribution of the variables $X_1, \ldots, X_p$ conditional on $Z = i$ is then given by a $p$-dimensional $\kappa_1 \times \cdots \times \kappa_p$ table

$$\mathbb{P}_i = \bigotimes_{j=1}^{p} \mathbf{p}_{ij},$$

whose $(l_1, l_2, \ldots, l_p)$-entry is $\prod_{j=1}^{p} \mathbf{p}_{ij}(l_j)$. Let

$$(1) \qquad \mathbb{P} = \sum_{i=1}^{r} \pi_i \mathbb{P}_i.$$

Then $\mathbb{P}$ is the distribution of a finite mixture of finite measure products, with a known number $r$ of components. The $\pi_i$ are interpreted as probabilities that a draw from the population is in the $i$th of $r$ classes. Conditioned on the class, the $p$ observable variables are independent. However, since the class is not discernible, the $p$ feature variables $X_j$ described by one-dimensional marginalizations of $\mathbb{P}$ are generally not independent.

We refer to the model described above as the $r$-class, $p$-feature model with state space $\{1, \ldots, \kappa_1\} \times \cdots \times \{1, \ldots, \kappa_p\}$, and denote it by $\mathcal{M}(r; \kappa_1, \kappa_2, \ldots, \kappa_p)$. Identifying the parameter space of this model with a subset $\Theta$ of $[0,1]^L$ where $L = (r-1) + r \sum_{i=1}^{p} (\kappa_i - 1)$ and letting $K = \prod_{i=1}^{p} \kappa_i$, the parameterization map for this model is

$$\Psi_{r,p,(\kappa_i)} : \Theta \to [0,1]^K.$$

In the following, we specify parameters by vectors such as $\boldsymbol{\pi}$ and $\mathbf{p}_{ij}$, always implicitly assuming their entries add to 1.

As previously noted, this model's parameters are not strictly identifiable if $r > 1$, since the sum in (1) can always be reordered without changing $\mathbb{P}$. Even modulo this label swapping, there are certainly special instances when identifiability will not hold. For instance, if $\mathbb{P}_i = \mathbb{P}_j$, then the parameters $\pi_i$ and $\pi_j$ can be varied, as long as their sum $\pi_i + \pi_j$ is held fixed, without effect on the distribution $\mathbb{P}$. Slightly more elaborate "special" instances of nonidentifiability can be constructed, but in full generality, this issue remains poorly understood. Ideally, one would know for which choices of $r, p, (\kappa_i)$, generic values of the model's parameters are identifiable up to permutation of the terms in (1), and, additionally, have a characterization of the exceptional set of parameters on which identifiability fails.

## 4. Kruskal's theorem and its consequences.
The basic identifiability result on which we build our later arguments is a result of Kruskal [29, 30] in the context of factor analyses for $p = 3$ features. Kruskal's result deals with a 3-way contingency table (or array) which cross-classifies a sample of $n$ individuals with respect to 3 polytomous variables (the $i$th of which takes values in $\{1, \ldots, \kappa_i\}$). If there is some latent variable $Z$ with values in $\{1, \ldots, r\}$ so that each of the $n$ individuals belongs to one of the $r$ latent classes and within the $l$th latent class, the 3 observed variables are mutually independent, then this $r$-class latent structure would serve as a simple explanation of the observed relationships among the variables in the 3-way contingency table. This latent structure analysis corresponds exactly to the model $\mathcal{M}(r; \kappa_1, \kappa_2, \kappa_3)$ described in the previous section.

To emphasize the focus on 3-variate models, note that in [30] Kruskal points out that 2-way tables arising from the model $\mathcal{M}(r; \kappa_1, \kappa_2)$ do not have a unique decomposition when $r \geq 2$. This nonidentifiability is intimately related to the nonuniqueness of certain matrix factorizations. While Goodman [22] studied the model $\mathcal{M}(r; \kappa_1, \kappa_2, \kappa_3, \kappa_4)$ for fitting to 4-way contingency tables, no formal result about uniqueness of the decomposition was established. In fact, nonidentifiability of the model under certain circumstances is highlighted in that work.

To present Kruskal's result, we introduce some algebraic notation. For $j = 1, 2, 3$, let $M_j$ be a matrix of size $r \times \kappa_j$, with $\mathbf{m}_i^j = (m_i^j(1), \ldots, m_i^j(\kappa_j))$ the $i$th row of $M_j$. Let $[M_1, M_2, M_3]$ denote the $\kappa_1 \times \kappa_2 \times \kappa_3$ tensor defined by

$$[M_1, M_2, M_3] = \sum_{i=1}^{r} \mathbf{m}_i^1 \otimes \mathbf{m}_i^2 \otimes \mathbf{m}_i^3.$$

In other words, $[M_1, M_2, M_3]$ is a three-dimensional array whose $(u, v, w)$ element is

$$[M_1, M_2, M_3]_{u,v,w} = \sum_{i=1}^{r} m_i^1(u) m_i^2(v) m_i^3(w)$$

for any $1 \leq u \leq \kappa_1, 1 \leq v \leq \kappa_2, 1 \leq w \leq \kappa_3$. Note that $[M_1, M_2, M_3]$ is left unchanged by simultaneously permuting the rows of all the $M_j$ and/or rescaling the rows so that the product of the scaling factors used for the $\mathbf{m}_i^j$, $j = 1, 2, 3$, is equal to 1.

A key point is that the probability distribution in a finite latent-class model with three observed variables is exactly described by such a tensor: let $M_j$, $j = 1, 2, 3$, be the matrix whose $i$th row is $\mathbf{p}_{ij} = \mathbb{P}(X_j = \cdot \mid Z = i)$. Let $\tilde{M}_1 = \text{diag}(\boldsymbol{\pi}) M_1$ be the matrix whose $i$th row is $\pi_i \mathbf{p}_{i1}$. Then the $(u, v, w)$ element of the tensor $[\tilde{M}_1, M_2, M_3]$ equals $\mathbb{P}(X_1 = u, X_2 = v, X_3 = w)$. Thus knowledge of the distribution of $(X_1, X_2, X_3)$ is equivalent to the knowledge of the tensor $[\tilde{M}_1, M_2, M_3]$. Note that the $M_i$'s are stochastic matrices, and thus the vector of $\pi_i$'s can be thought of as scaling factors.

For a matrix $M$, the *Kruskal rank* of $M$ will mean the largest number $I$ such that every set of $I$ rows of $M$ are independent. Note that this concept would change if we replaced "row" by "column," but we will only use the row version in this paper. With the Kruskal rank of $M$ denoted by $\text{rank}_K M$, we have

$$\text{rank}_K M \leq \text{rank} M$$

and equality of rank and Kruskal rank does not hold in general. However, in the particular case where a matrix $M$ of size $p \times q$ has rank $p$, it also has Kruskal rank $p$.

The fundamental algebraic result of Kruskal is the following.

THEOREM 1 (Kruskal [29, 30]). *Let* $I_j = \text{rank}_K M_j$. *If*

$$I_1 + I_2 + I_3 \geq 2r + 2,$$

*then* $[M_1, M_2, M_3]$ *uniquely determines the* $M_j$, *up to simultaneous permutation and rescaling of the rows.*

The equivalence between the distributions of 3-variate latent-class models and 3-tensors, combined with the fact that rows of stochastic matrices sum to 1, gives the following reformulation.

COROLLARY 2. *Consider the model* $\mathcal{M}(r; \kappa_1, \kappa_2, \kappa_3)$, *with the parameterization of Section 3. Suppose all entries of* $\boldsymbol{\pi}$ *are positive. For each* $j = 1, 2, 3$, *let* $M_j$ *denote the matrix whose rows are* $\mathbf{p}_{ij}$, $i = 1, \ldots, r$, *and let* $I_j$ *denote its Kruskal rank. Then if*

$$I_1 + I_2 + I_3 \geq 2r + 2,$$

*the parameters of the model are uniquely identifiable, up to label swapping.*

By observing that Kruskal's condition on the sum of Kruskal ranks can be expressed through polynomial inequalities in the parameters, and thus holds generically, we obtain the following corollary.

COROLLARY 3. *The parameters of the model* $\mathcal{M}(r; \kappa_1, \kappa_2, \kappa_3)$ *are generically identifiable, up to label swapping, provided*

$$\min(r, \kappa_1) + \min(r, \kappa_2) + \min(r, \kappa_3) \geq 2r + 2.$$

*The assertion remains valid if, in addition, the class proportions* $\{\pi_i\}_{1 \leq i \leq r}$ *are held fixed and positive in the model.*

For the last statement of this corollary, we note that if the mixing proportions are positive then one can translate Kruskal's condition into a polynomial requirement so that only the parameters $\mathbf{p}_{ij}(l) = \mathbb{P}(X_j = l \mid Z = i)$ appear. Thus the generic aspect only concerns this part of the parameter space, and not the part with the proportions $\pi_i$. As a consequence, the statement is valid when the proportions are held fixed in $(0, 1)$. This is of great importance, as often statisticians assume that these proportions are fixed and known (for instance using $\pi_i = 1/r$ for every $1 \leq i \leq r$). Without observing this fact, we would not have a useful identifiability result in the case of known $\pi_i$, since fixing values of the $\pi_i$ results in considering a subvariety of the full parameter space, which a priori might be included in the subvariety of nonidentifiable parameters allowed by Corollary 3.

**5. Parameter identifiability of finite mixtures of finite measure products.** Finite mixtures of products of finite measure are widely used to model data, for instance in biological taxonomy, medical diagnosis or classification of text documents [21, 35]. The identifiability issue for these models was first addressed forty years ago by Teicher [42]. Teicher's result states the equivalence between identifiability of mixtures of product measure distributions and identifiability of the corresponding one-dimensional mixture models. As a consequence, finite mixtures of Bernoulli products are not identifiable in a strict sense [23]. Teicher's result is valid for finite mixtures with an unknown number of components, but it can easily be seen that nonidentifiability occurs even with a known number of components [6], Section 1. The very simplicity of the equivalence condition stated by Teicher [42] likely impeded statisticians from looking further at this issue.

Here we prove that finite mixtures of Bernoulli products (with a known number of components) are in fact generically identifiable, indicating why these models are well behaved in practice with respect to statistical parameter inference, despite their lack of strict identifiability [6].

To obtain our results, we must first pass from Kruskal's theorem on 3-variate models to a similar one for $p$-variate models. To do this, we observe that $p$ observed variables can be combined into 3 agglomerate variables, so that Kruskal's result can be applied.

THEOREM 4. *Consider the model* $\mathcal{M}(r; k_1, \ldots, k_p)$ *where* $p \geq 3$. *Suppose there exists a tripartition of the set* $S = \{1, \ldots, p\}$ *into three disjoint nonempty subsets* $S_1, S_2, S_3$, *such that if* $\kappa_i = \prod_{j \in S_i} k_j$ *then*

$$(2) \qquad \min(r, \kappa_1) + \min(r, \kappa_2) + \min(r, \kappa_3) \geq 2r + 2.$$

*Then model parameters are generically identifiable, up to label swapping. Moreover, the statement remains valid when the mixing proportions* $\{\pi_i\}_{1 \leq i \leq r}$ *are held fixed and positive.*

Considering the special case of finite mixtures of $r$ Bernoulli products with $p$ components [i.e., the $r$-class, $p$-binary feature model $\mathcal{M}(r; 2, 2, \ldots, 2)$], to obtain the strongest identifiability result, we choose a tripartition that maximizes the left-hand side of inequality (2). Doing so yields the following.

COROLLARY 5. *Parameters of the finite mixture of* $r$ *different Bernoulli products with* $p$ *components are generically identifiable, up to label swapping, provided*

$$p \geq 2\lceil \log_2 r \rceil + 1,$$

*where* $\lceil x \rceil$ *is the smallest integer at least as large as* $x$.

Note that generic identifiability of this model for sufficiently large values of $p$ is a consequence of the results of Elmore, Hall and Neeman, in [14], although neither the generic nature of the result, nor the fact that the model is simply a mixture of Bernoulli products, is noted by the authors. Moreover, our lower bound on $p$ to ensure generic identifiability is superior to the one obtained in [14]. Indeed, letting $C(r)$ be the minimal integer such that if $p \geq C(r)$ then the $r$-class, $p$-binary feature model is generically identifiable, then [14] established that

$$\log_2 r \leq C(r) \leq c_2 r \log_2 r$$

for some effectively computable constant $c_2$. While the lower bound for $C(r)$ is easy to obtain from the necessity that the dimension of the parameter space, $rp + (r-1)$, be no larger than that of the distribution space $2^p - 1$, the upper bound required substantial work. Corollary 5 above establishes the stronger result that

$$C(r) \leq 2\lceil \log_2 r \rceil + 1.$$

Note that this new upper bound, along with the simple lower bound, shows that the order of growth of $C(r)$ is precisely $\log_2 r$.

For the more general $\mathcal{M}(r; \kappa, \ldots, \kappa)$ model, our lower bound on the number of variates needed to generically identify the parameters, up to label swapping, is

$$p \geq 2\lceil \log_\kappa r \rceil + 1.$$

The proof of this bound follows the same lines as that of Corollary 5, and is therefore omitted.

**6. Hidden classes models with dependent observations.** In this section, we give several additional illustrations of the applicability of Kruskal's result in the context of dependent observations. The hidden Markov models and random graph mixture models we consider may at first appear to be far from the focus of Kruskal's theorem. This is not the case, however, as in both the observable variables are independent when appropriately conditioned on hidden ones. We succeed in embedding these models into an appropriate $\mathcal{M}(r; \kappa_1, \kappa_2, \kappa_3)$ and then use extensive algebraic arguments to obtain the (generic) identifiability results we desire.

6.1. *Hidden Markov models.* Almost 40 years ago, Petrie ([39], Theorem 1.3) proved generic identifiability, up to label swapping, for discrete hidden Markov models (HMMs). We offer a new proof, based on Kruskal's theorem, of this well-known result. This provides an interesting alternative to Petrie's more direct approach, and one that might extend to more complex frameworks, such as Markov chains with Markov regime, where no identifiability results are known (see, for instance, [9]). Moreover, as a by-product, our approach establishes a new bound on the number of consecutive variables needed, such that the marginal distribution for a generic HMM uniquely determines the full probability distribution.

We first briefly describe HMMs. Consider a stationary Markov chain $\{Z_n\}_{n \geq 0}$ on state space $\{1, \ldots, r\}$ with transition matrix $A$ and initial distribution $\boldsymbol{\pi}$ (assumed to be the stationary distribution). Conditional on $\{Z_n\}_{n \geq 0}$, the observations $\{X_n\}_{n \geq 0}$ on state space $\{1, \ldots, \kappa\}$ are assumed to be i.i.d., and the distribution of each $X_n$ only depends on $Z_n$. Denote by $B$ the matrix of size $r \times \kappa$ containing the conditional probabilities $\mathbb{P}(X_n = k \mid Z_n = i)$. The process $\{X_n\}_{n \geq 0}$ is then a hidden Markov chain. Note that this is not a Markov process. The matrices $A$ and $B$ constitute the parameters for the $r$-hidden state, $\kappa$-observable state HMM, and the parameter space can be identified with a full-dimensional subset of $\mathbb{R}^{r(r + \kappa - 2)}$. We refer to [5, 16] for more details on HMMs.

Petrie [39] describes quite explicitly, for fixed $r$ and $\kappa$, a subvariety of the parameter space for an HMM on which identifiability might fail. Indeed, Petrie proved that the set of parameters on which identifiability holds is the intersection of the following: the set of regular HMMs; the set where the components of the matrix $B$, namely $\mathbb{P}(X_n = k \mid Z_n = i)$ are nonzero; the set where some row of $B$ has distinct entries [namely there exists some $i \in \{1, \ldots, r\}$ such that all the $\{\mathbb{P}(X_n = k \mid Z_n = i)\}_k$ are distinct]; the set where the matrix $A$ is nonsingular, and 1 is an eigenvalue with multiplicity one for $A$ [namely, $P'(1, A) \neq 0$ where $P(\lambda, A) = \det(\lambda I - A)$]. Regular HMMs were first described by Gilbert in [20]. The definition relies on a notion of *rank* and an HMM is regular if its rank is equal to its number of hidden states $r$. More details may be found in [17, 20].

The result of Petrie assumes knowledge of the whole probability distribution of the HMM. But it is known ([17], Lemma 1.2.4) that the distribution of an HMM with $r$ hidden states and $\kappa$ observed states, is completely determined by the marginal distribution of $2r$ consecutive variables. An even stronger result appears in [38], Chapter 1, Corollary 3.4: the marginal distribution of $2r - 1$ consecutive variables suffices to reconstruct the whole HMM distribution. Combining these results shows that generic identifiability holds for HMMs from the distribution of $2r - 1$ consecutive observations. Note there is no dependence of this number on $\kappa$, even though one might suspect a larger observable state space would aid identifiability.

Using Kruskal's theorem we prove the following.

THEOREM 6. *The parameters of an HMM with $r$ hidden states and $\kappa$ observable states are generically identifiable from the marginal distribution of $2k + 1$ consecutive variables provided $k$ satisfies*

$$(3) \qquad \binom{k + \kappa - 1}{\kappa - 1} \geq r.$$

While we do not explicitly characterize a set of possibly nonidentifiable parameters as Petrie did, in principle we could do so.

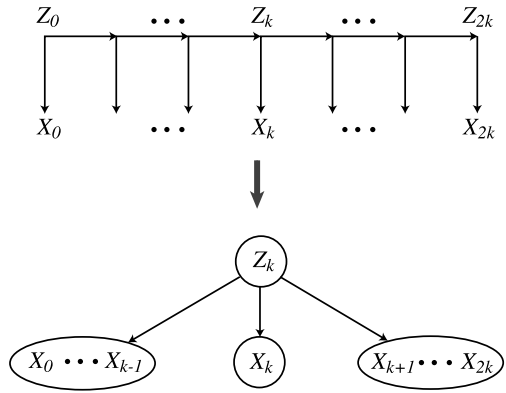

Fig. 1. *Embedding the hidden Markov model into a simpler latent-class model.*

Note, however, that we require only the marginal of $2k+1$ consecutive variables, where $k$ satisfies an explicit condition involving $\kappa$. The worst case (i.e., the largest value for $k$) arises when $\kappa = 2$, since $\kappa \mapsto \binom{k+\kappa-1}{\kappa-1}$ is an increasing function for positive $\kappa$. In this worst case, we easily compute that $2k+1 = 2r-1$ consecutive variables suffice to generically identify the parameters. Thus our approach yields generic versions of the claims of [17] and [38] described above.

Moreover, when the number $\kappa$ of observed states is increased, the minimal value of $2k+1$ which ensures identifiability by Theorem 6 becomes smaller. Thus the fact that generic HMMs are characterized by the marginal of $2k+1$ consecutive variables, where $k$ satisfies (3), results in a much better bound than $2r-1$ as soon as the observed state space has more than 2 points. In this sense, our result is stronger than the one of Paz [38].

In proving Theorem 6, we embed the hidden Markov model in a simpler latent-class model, as illustrated in Figure 1. The hidden variable $Z_k$ is the only one we preserve, while we cluster the observed variables into groups so they may be treated as 3 observed variables. According to properties of graphical models (see, e.g., [31]), the agglomerated variables are independent when conditioned on $Z_k$. So Kruskal's theorem applies. However, additional algebraic arguments are needed to see that the embedding gives sufficiently generic points so that we may identify parameters.

To conclude this section, note that Leroux [32] used Teicher's result [42] to establish a sufficient condition for exact identifiability of parametric hidden Markov models, with possibly continuous observations.

6.2. *A random graph mixture model.* We next illustrate the application of our method to studying a random graph mixture model. This heterogeneous model is used in a wide range of applications, such as molecular biology (gene interactions or metabolic networks [12]), social sciences (relationships or co-authorship networks [36]) or the study of the world wide web (hyperlinks graphs [44]).

We consider a random graph mixture model in which each node belongs to some unobserved class, or group, and conditional on the classes of all nodes, the edges are independent random variables whose distributions depend only on the classes of the nodes they connect. More precisely, consider an undirected graph with $n$ nodes labeled $1, \ldots, n$ and where presence/absence of an edge between two different nodes $i$ and $j$ is given by the indicator variable $X_{ij}$. Let $\{Z_i\}_{1 \leq i \leq n}$ be i.i.d. random variables with values in $\{1, \ldots, r\}$ and probability distribution $\boldsymbol{\pi} \in (0, 1)^r$ representing node classes. Conditional on the classes of the nodes $\{Z_i\}$, the edge indicators $X_{ij}$ are independent random variables whose distribution is Bernoulli with some parameter $p_{Z_i Z_j}$. The between-groups connection parameters $p_{ij} \in [0, 1]$ satisfy $p_{ij} = p_{ji}$, for all $1 \leq i, j \leq r$. We emphasize that the observed random variables $X_{ij}$ for this model are not independent, just as the observed variables were not independent in the mixture models considered earlier in this paper.

The interest in the random graph model lies in the fact that different nodes may have different connectivity properties. For instance one class may describe *hubs* which are nodes with a very high connectivity, and a second class may contain the others nodes with a lower overall connectivity. Thus one can model different node behaviours with a reasonable number of parameters. Examples of networks easily modelled with this approach, and more details on the properties of this model, can be found in [12].

This model has been rediscovered many times in the literature and in various fields of applications. A nonexhaustive bibliography includes [12, 18, 36, 41]. However, identifiability of the parameters for this random graph model has never been addressed in the literature.

Frank and Harary [18] study the statistical inference of the parameters in the restricted $\alpha$–$\beta$ or affiliation model. In this setup, only two parameters are used to model the intra-group and inter-group probabilities of an edge occurrence $p_{ii} = \alpha$, $1 \leq i \leq r$ and $p_{ij} = \beta$, $1 \leq i < j \leq r$. Using the total number of edges, the proportion of transitive triads and the proportion of 3-cycles among triads (see definitions (3) and (4) in [18]), they obtain estimates of the parameters $\alpha, \beta$ and sometimes $r$, in various cases ($\alpha, \beta$ unknown and $\pi_i = 1/r$ with unknown $r$, for instance). However, they do not discuss the uniqueness of the solutions of the nonlinear equations defining those estimates (see (16), (29) and (33) in [18]).

We prove the following result.

THEOREM 7. *The parameters of the random graph model with $r = 2$ node states are strictly identifiable, up to label swapping, provided there are at least 16 nodes and the connection parameters $\{p_{11}, p_{12}, p_{22}\}$ are distinct.*

Our basic approach is to embed the random graph model into a model to which Kruskal's theorem applies. Since we have many hidden variables, one for each node, we combine them into a single composite variable describing the states of all nodes at once. Since the observed edge variables are binary, we also combine them into collections, to create 3 composite edge variables with more states. We do this in such a way that the composite edge variables are still independent conditioned on the composite node state. The main technical difficulty is that the matrices whose entries give probabilities of observing a composite edge variable conditioned on a composite node state must have well-understood Kruskal rank. This requires some involved algebraic work.

The random graph model will be studied more thoroughly in a forthcoming work. The special case of the affiliation model, which is not adressed by Theorem 7, will be dealt with there as well.

**7. Finite mixtures of nonparametric measure products.** In this section, we consider a nonparametric model of finite mixtures of $r$ different probability distributions $\mu_1, \ldots, \mu_r$ on $\mathbb{R}^p$, with $p \geq 3$. For every $1 \leq i \leq r$, we denote by $\mu_i^j$ the $j$th marginal of $\mu_i$ and $F_i^j$ the corresponding c.d.f. (defined by $F_i^j(t) = \mu_i^j((-\infty, t])$ for any $t \in \mathbb{R}$). Without loss of generality, we may assume that the functions $F_i^j$ are absolutely continuous.

For our first result, we assume further that the mixture model has the form

$$(4) \qquad \mathbb{P} = \sum_{i=1}^r \pi_i \mu_i = \sum_{i=1}^r \pi_i \prod_{j=1}^p \mu_i^j,$$

in which, conditional on a latent structure (specified by the proportions $\pi_i$), the $p$ variates are independent. The $\mu_i^j$ are viewed in a nonparametric setting.

In the next theorem, we prove identifiability of the model's parameters—that is, that $\mathbb{P}$ uniquely determines the factors appearing in (4)—under a mild and explicit regularity condition on $\mathbb{P}$, as soon as there are at least 3 variates and $r$ is known. Making a judicious use of cut points to discretize the distribution, and then using Kruskal's work, we prove the following.

THEOREM 8. *Let $\mathbb{P}$ be a mixture of the form (4), such that for every $j \in \{1, \ldots, p\}$, the measures $\{\mu_i^j\}_{1 \leq i \leq r}$ are linearly independent. Then, if $p \geq 3$, the parameters $\{\pi_i, \mu_i^j\}_{1 \leq i \leq r, 1 \leq j \leq p}$ are strictly identifiable from $\mathbb{P}$, up to label swapping.*

This result also generalizes to nonparametric mixture models where at least three blocks of variates are independent conditioned on the latent structure. Let $b_1, \ldots, b_p$ be integers with $p \geq 3$ the number of blocks, and the $\mu_i^j$ be absolutely continuous probability measures on $\mathbb{R}^{b_j}$. With $m = \sum_j b_j$, consider the mixture distribution on $\mathbb{R}^m$ given by

$$\mathbb{P} = \sum_{i=1}^{r} \pi_i \prod_{j=1}^{p} \mu_i^j. \tag{5}$$

THEOREM 9. *Let $\mathbb{P}$ be a mixture of the form (5), such that for every $j \in \{1, \ldots, p\}$, the measures $\{\mu_i^j\}_{1 \leq i \leq r}$ on $\mathbb{R}^{b_j}$ are linearly independent. Then, if $p \geq 3$, the parameters $\{\pi_i, \mu_i^j\}_{1 \leq i \leq r, 1 \leq j \leq p}$ are strictly identifiable from $\mathbb{P}$, up to label swapping.*

Both Theorems 8 and 9 could be strengthened somewhat, as their proofs do not depend on the full power of Kruskal's theorem. As an analog of Kruskal rank for a matrix, say a finite set of measures has *Kruskal rank* $k$, if $k$ is the maximal integer such that every $k$-element subset is linearly independent. Then, for instance, when $p = 3$, straightforward modifications of the proofs establish identifiability provided the sum of the Kruskal ranks of the sets $\{\mu_i^j\}_{1 \leq i \leq r}$ for $j = 1, 2, 3$ is at least $2r + 2$.

Note that an earlier result linking identifiability with linear independence of the densities to be mixed appears in [43] in a parametric context. Combining this statement with the one obtained by Teicher [42], we get that in the parametric framework, a sufficient and necessary condition for strict identifiability of finite mixtures of product distributions is the linear independence of the univariate components (this statement does not require the knowledge of $r$). In this sense, our result may be seen as a generalization of these statements in the nonparametric context.

Our results should be compared to two previous ones. First, Hettmansperger and Thomas [26] studied the identical marginals case, where for each $1 \leq i \leq r$, we have $\mu_i^j = \mu_i^k$ for all $1 \leq j, k \leq p$ (with corresponding c.d.f. $F_i$). They proved that, as soon as $p \geq 2r - 1$, and there exists some $c \in \mathbb{R}$, such that the $\{F_i(c)\}_{1 \leq i \leq r}$ are all different, the mixing proportions $\pi_i$ are identifiable. Although they do not state it (because they are primarily interested in estimation procedures and not identifiability), they also identify the c.d.f.s $F_i(c), 1 \leq i \leq r$, at any point $c \in \mathbb{R}$ such that the $\{F_i(c)\}_{1 \leq i \leq r}$ are all different.

Second, Hall and Zhou [25] proved that for $r = 2$, if $p \geq 3$ and the mixture $\mathbb{P}$ is such that its two-dimensional marginals are not the product of the corresponding one-dimensional ones, that is, that for any $1 \leq j, k \leq p$ we have

$$\mathbb{P}(X_j, X_k) \neq \mathbb{P}(X_j)\mathbb{P}(X_k), \tag{6}$$

then the parameters $\boldsymbol{\pi}$ and $F_i^j$ are uniquely identified by $\mathbb{P}$. They also provide consistent estimation procedures for the mixing proportions as well as for the univariate c.d.f.'s $\{F_i^j\}_{1 \leq j \leq p, i=1,2}$.

Hall and Zhou state that their results, "apparently do not have straightforward generalizations to mixtures of three or more products of independent components." In fact, Theorem 8 can already be viewed as such a generalization, since in the $r = 2$ case we will show that inequality (6) is in fact equivalent to the independence for each $j$ of the set $\{\mu_i^j\}_{1 \leq i \leq 2}$.

To develop a more direct generalization of the condition of Hall and Zhou, we say a bivariate continuous probability distribution is of *rank $r$* if it can be written as a sum of $r$ products of signed univariate distributions (not necessarily probability distributions), but no fewer. We emphasize that we allow the univariate distributions to have negative values, even though the bivariate does not. The related notion of *nonnegative rank*, which additionally requires that the univariate distributions be nonnegative, will not play a role here.

This definition of rank of a bivariate distribution is a direct generalization of the notion of the rank for a matrix, with the bivariate distribution replacing the matrix, the univariate ones replacing vectors, and the product of distributions replacing the vector outer product. Moreover, in the case $r = 2$, inequality (6) is equivalent to saying $\mathbb{P}(X_j, X_k)$ has rank 2, since its rank is at most 2 from the expression (4), and if it had rank 1 then marginalizing would show $\mathbb{P}(X_j, X_k) = \mathbb{P}(X_j)\mathbb{P}(X_k)$.

Next we connect this concept to the hypotheses of Theorem 8.

LEMMA 10. *Consider a bivariate distribution of the form*

$$\mathbb{P}(X_1, X_2) = \sum_{i=1}^r \pi_i \mu_i^1(X_1) \mu_i^2(X_2).$$

*Then $\mathbb{P}(X_1, X_2)$ has rank $r$ if, and only if, for each of $j = 1, 2$ the measures $\{\mu_i^j\}_{1 \leq i \leq r}$ are linearly independent.*

From Theorem 8, this immediately yields the following.

COROLLARY 11. *Let $\mathbb{P}$ be a mixture of the form (4), and suppose that for every $j \in \{1, \ldots, p\}$, there is some $k \in \{1, \ldots, p\}$ such that the marginal $\mathbb{P}(X_j, X_k)$ is of rank $r$. Then, if $p \geq 3$, the parameters $\{\pi_i, \mu_i^j\}_{1 \leq i \leq r, 1 \leq j \leq p}$ are strictly identifiable from $\mathbb{P}$, up to label swapping.*

Note that if the distribution $\mathbb{P}$ arising from (4) could be written with strictly fewer than $r$ different product distributions, then the assumption of Corollary 11 (and of Theorem 8) would not be met. Here, we state that a slightly stronger condition, namely that any two-variate marginal of $\mathbb{P}$ cannot be written as the sum of strictly fewer than $r$ product components, suffices to ensure identifiability of the parameters. Note also that the condition appearing in Theorem 8 is stated in terms of the parameters of the distribution $\mathbb{P}$, which are unknown to the statistician. However, the rephrased assumption appearing in Corollary 11 is stated in terms of the observation distribution. Thus one could imagine testing this assumption on the data (even if this might be a tough issue) prior to statistical inference.

The restriction to $p \geq 3$, which arises in our method from using Kruskal's theorem, is necessary in this context. Indeed, in the case of $r = 2$ groups and $p = 2$ variates, [25] proved that there exists a two-parameter continuum of points $(\boldsymbol{\pi}, \{\mu_i^j\}_{1 \leq i \leq 2, 1 \leq j \leq 2})$ solving equation (4). This simply echos the nonidentifiability in the case of 2 variates with finite state spaces commented on in [30].

For models with more than 2 components in the mixture, Benaglia, Chauveau and Hunter [4] recently proposed an algorithmic estimation procedure, without insurance that the model would be identifiable. Our results states that under mild regularity conditions, it is possible to identify the parameters from the mixture $\mathbb{P}$, at least when the number of components $r$ is fixed. Thus our approach gives some theoretical support to the procedure developed in [4].

Finally, recall that in [14, 24], an upper bound on the number of variates needed to ensure "generic" identifiability of the model is claimed which is of the order $r \log_2(r)$. Our Theorem 8 lowers this bound considerably, as it shows 3 variates suffice to identify the model, regardless of the value of $r$ (under a mild regularity assumption).

## 8. Proofs.

PROOF OF COROLLARY 2. For each $j = 1, 2, 3$, let $M_j$ be the matrix of size $r \times \kappa_j$ describing the probability distribution of $X_j$ conditional on $Z$. More precisely, the $i$th row of $M_j$ is $\mathbf{p}_{ij} = \mathbb{P}(X_j = \cdot \mid Z = i)$, for any $i \in \{1, \ldots, r\}$. Let $\tilde{M}_1$ be the matrix of size $r \times \kappa_1$ such that its $i$th row is $\pi_i \mathbf{p}_{i1} = \pi_i \mathbb{P}(X_1 = \cdot \mid Z = i)$. Kruskal ranks of $M_j$ and $\tilde{M}_1$ are denoted $I_j$ and $\tilde{I}_1$, respectively. We have already seen that the tensor $[\tilde{M}_1, M_2, M_3]$ describes the probability distribution of the observations $(X_1, X_2, X_3)$. Kruskal's result states that, as soon as Kruskal ranks satisfy the condition $\tilde{I}_1 + I_2 + I_3 \geq 2r + 2$, this probability distribution uniquely determines the matrices $\tilde{M}_1, M_2$ and $M_3$ up to rescaling of the rows and label swapping. Note that as $\boldsymbol{\pi}$ has positive entries, Kruskal rank $\tilde{I}_1$ is equal to $I_1$. Moreover, using that the matrices $M_1, M_2$ and $M_3$ are stochastic and that the entries of $\boldsymbol{\pi}$ are positive, the corollary follows. □

PROOF OF COROLLARY 3.   We first show that for any fixed choice of
a positive integer $I_j \leq \min(r, \kappa_j)$, those $r \times \kappa_j$ matrices $M_j$, whose Kruskal
rank is strictly less than $I_j$, form a proper algebraic variety. But the matrices
for which a specific set of $I_j$ rows are dependent is the zero set of all $I_j \times I_j$
minors obtained from those rows. By taking appropriate products of these
minors for all such sets of $I_j$ rows we may construct a set of polynomials
whose zero set is precisely those matrices of Kruskal rank less than $I_j$. This
variety is proper, since matrices of full rank do not lie in it.

Thus the set of triples of matrices $(M_1, M_2, M_3)$ for which the Kruskal
rank of $M_i$ is strictly less than $\min(r, \kappa_i)$ forms a proper subvariety. For
triples not in this subvariety, our assumptions ensure that the rank inequal-
ity of Corollary 3 holds, so the inequality holds generically. If the $\pi_i$'s are
fixed and positive, the proof is complete. Otherwise, note that the set of
parameters with vectors $\boldsymbol{\pi}$ admitting zero entries is also a proper subvariety
of the parameter set.   □

PROOF OF THEOREM 4.   Our goal is to apply Kruskal's result to models
with more than 3 observed variables by means of a "grouping" argument.
We require a series of lemmas to accomplish this.

First, given an $n \times a_1$ matrix $A_1$ and an $n \times a_2$ matrix $A_2$, define the
$n \times a_1 a_2$ matrix $A = A_1 \otimes^{\text{row}} A_2$, as the row-wise tensor product, so that

$$A(i, a_2(j-1) + k) = A_1(i, j) A_2(i, k).$$

The proof of the following lemma is straightforward and therefore omitted.

LEMMA 12.   *If conditional on a finite random variable $Z$, the random
variables $X_1, X_2$ are independent, with the distribution of $X_i$ conditional
on $Z$ given by the matrix $A_i$ of size $r \times a_i$, then the row tensor product
$A = A_1 \otimes^{\text{row}} A_2$ of size $r \times (a_1 a_2)$ contains the probability distribution of
$(X_1, X_2)$ conditional on $Z$.*

For each $j \in \{1, \ldots, p\}$, denote by $M_j$ the $r \times k_j$ matrix whose $i$th row
is $\mathbb{P}(X_j = \cdot \mid Z = i)$. Introduce three matrices $N_i$, $i = 1, 2, 3$, of size $r \times \kappa_i$,
defined as

$$N_i = \bigotimes_{j \in S_i}^{\text{row}} M_j$$

and the tensor $N = [\tilde{N}_1, N_2, N_3]$, where the $i$th row of $\tilde{N}_1$ is $\pi_i$ times the
$i$th row of $N_1$. According to Lemma 12, the tensor $N$ contains the proba-
bilities of the three clumped variables $(\{X_j\}_{j \in S_1}, \{X_j\}_{j \in S_2}, \{X_j\}_{j \in S_3})$. Thus
knowledge of the distribution of the observations is equivalent to knowledge
of $N$. Moreover, for parameters $\boldsymbol{\pi}$ having positive entries (which is a generic
condition), the Kruskal ranks of $\tilde{N}_1$ and $N_1$ are equal.

In the next lemma we characterize the Kruskal rank of the row-tensor product obtained from generic matrices $A_i$.

LEMMA 13.   *Let $A_i$, $i = 1, \ldots, q$, denote $r \times a_i$ matrices, $a = \prod_{i=1}^{q} a_i$ and*

$$A = \bigotimes_{i=1,\ldots,q}^{\mathrm{row}} A_i,$$

*the $r \times a$ matrix obtained by taking tensor products of the corresponding rows of the $A_i$. Then for generic $A_i$'s,*

$$\mathrm{rank}_K A = \mathrm{rank}\, A = \min(r, a).$$

PROOF.   The condition that a matrix $A$ not have full rank (resp. full Kruskal rank) is equivalent to the simultaneous vanishing of its maximal minors (resp. *idem* when $r \leq a$, and equivalent to the existence of one vanishing maximal minor when $r > a$). Composing the map sending $\{A_i\} \to A$ with these minors gives polynomials in the entries of the $A_i$. To see that the polynomials in the entries of the $A_i$ are nonzero, it is enough to exhibit a single choice of the $A_i$ for which $A$ has full rank (resp. full Kruskal rank).

Let $x_{ij}$, $i = 1, \ldots, q$, $j = 1, \ldots, a_i$, be distinct prime numbers. Consider $A_i$ defined by

$$A_i = \begin{pmatrix} 1 & 1 & \cdots & 1 \\ x_{i1} & x_{i2} & \cdots & x_{ia_i} \\ x_{i1}^2 & x_{i2}^2 & \cdots & x_{ia_i}^2 \\ \vdots & \vdots & \ddots & \vdots \\ x_{i1}^{r-1} & x_{i2}^{r-1} & \cdots & x_{ia_i}^{r-1} \end{pmatrix}.$$

For any vector $\mathbf{y} \in \mathbb{C}^t$, let $W(\mathbf{y}) = W(y_1, y_2, \ldots, y_t)$ denote the $t \times t$ Vandermonde matrix, with entries $y_j^{i-1}$.

Suppose first that $a \geq r$. Then the rows of $A$ are the first $r$ rows of the Vandermonde matrix $W(\tilde{\mathbf{y}})$, where $\tilde{\mathbf{y}}$ is a vector whose entries are $\prod_i x_{ij_i}$ for choices of $1 \leq j_i \leq a_i$. As the products $\prod_i x_{ij_i}$ are distinct by choice of the $x_{ij}$, $W(\tilde{\mathbf{y}})$ is nonsingular, so $A$ has rank and Kruskal rank equal to $r$.

If instead $r > a$, then the first $a$ rows of $A$ form an invertible Vandermonde matrix. Thus $A$ is of rank $a$. To argue that $A$ has full Kruskal rank, compose the map $\{A_i\} \to A$ with the $a \times a$ minor from the first $a$ rows of $A$. This gives us a polynomial in the entries of the $A_i$, the nonvanishing of which ensures the first $a$ rows of $A$ are independent. This polynomial is not identically zero, since a specific choice of the $A_i$'s such that the first $a$ rows of $A$ are independent has been given. By composing this polynomial with maps that permute rows of all the $A_i$ simultaneously, we may construct nonzero polynomials whose nonvanishing ensures all other sets of $a$ rows of $A$ are independent. The proper subvariety defined by the product of these polynomials, then, is precisely those choices of $\{A_i\}$ for which $A$ is not of full Kruskal rank. This concludes the proof of the lemma.   $\square$

Returning to the proof of Theorem 4, note that to apply the preceding lemma to stochastic matrices $M_j$, we must address the fact that each row of each $M_j$ sums to 1. However, as both rank and Kruskal rank are unaffected by multiplying rows by nonzero scalars, and rows sums being nonzero is a generic condition (defined by the nonvanishing of linear polynomials), we see immediately the conclusion of Lemma 13 holds when all the $M_j$ are additionally assumed to have row sums equal to 1.

We thus see that for generic $M_j$, the matrices $N_i$ defined above have Kruskal rank $I_i = \min(r, \kappa_i)$. Now by assumption, the matrices $N_i$ satisfy the condition of Corollary 2. This implies that the tensor $N = [\tilde{N}_1, N_2, N_3]$ uniquely determines the matrices $N_i$ and the vector $\boldsymbol{\pi}$, up to permutation of the rows. We need a last lemma before completing the proof of Theorem 4.

LEMMA 14.   *Suppose $A = \bigotimes_{i=1,\dots,q}^{\mathrm{row}} A_i$ where the $A_i$ are stochastic matrices. Then the $A_i$ are uniquely determined by $A$.*

PROOF.   Since each row of each $A_i$ sums to 1, one easily sees that each entry in $A_i$ can be recovered as a sum of certain entries in the same row of $A$.   $\square$

Using this lemma, we have that each $N_i$ uniquely determines the matrices $M_j$ for $j \in S_i$, and Theorem 4 follows.   $\square$

PROOF OF COROLLARY 5.   It is enough to consider the case where $p = 2\lceil \log_2 r \rceil + 1$. With $k = \lceil \log_2 r \rceil$, we have that $2^{k-1} < r \leq 2^k$. Choosing

$$\kappa_1 = \kappa_2 = 2^k, \qquad \kappa_3 = 2,$$

inequality (2) in Theorem 4 holds.   $\square$

PROOF OF THEOREM 6.   The $2k + 1$ consecutive observed variables can be taken to be $X_0, X_1, \dots, X_{2k}$. Note that the transition matrix from $Z_i$ to $Z_{i-1}$ is given by $A' = \mathrm{diag}(\boldsymbol{\pi})^{-1} A^T \mathrm{diag}(\boldsymbol{\pi})$.

Let $B_1$ be the $r \times \kappa^k$ matrix giving probabilities of joint states of $X_0, X_1, \dots, X_{k-1}$ conditioned on the states of $Z_k$. Similarly, let $B_2$ be the $r \times \kappa^k$ matrix giving probabilities of joint states of $X_{k+1}, \dots, X_{2k}$ conditioned on the states of $Z_k$.

Now the joint distribution for the model $\mathcal{M}(r; \kappa^k, \kappa^k, \kappa)$ with parameters $\boldsymbol{\pi}, B_1, B_2, B$ is the same as that of the HMM with parameters $\boldsymbol{\pi}, A, B$. Thus we apply Kruskal's theorem, after we first show $\boldsymbol{\pi}, B_1, B_2$ are sufficiently generic to do so for generic choices of $A, B$. The entries of $\boldsymbol{\pi}$ have been assumed to be positive. With $I_M$ denoting the Kruskal rank of a matrix $M$, in order to apply Corollary 2, we want to ensure

$$I_{B_1} + I_{B_2} + I_B \geq 2r + 2.$$

Making the generic assumption that $B$ has Kruskal rank at least 2, it is sufficient to make

$$I_{B_1}, I_{B_2} \geq r,$$

that is, to require that $B_1, B_2$ have full row rank.

Now $B_1, B_2$ can be explicitly given as

$$
\begin{aligned}
(7) \qquad B_1 &= A'(B \otimes^{\text{row}} (\cdots A'(B \otimes^{\text{row}} (A'(B \otimes^{\text{row}} (A'B)))) \cdots)), \\
B_2 &= A(B \otimes^{\text{row}} (\cdots A(B \otimes^{\text{row}} (A(B \otimes^{\text{row}} (AB)))) \cdots))
\end{aligned}
$$

with $k$ copies of $A'$ and of $B$ appearing in the expression for $B_1$, and $k$ copies of $A$ and $B$ appearing in that for $B_2$. To show these have full row rank for generic choices of stochastic $A$ and $B$, it is enough to show they have full row rank for some specific choice of stochastic $A, B, \boldsymbol{\pi}$, since that will establish that some $r \times r$ minors of $B_1, B_2$ are nonzero polynomials in the entries of $A', B$ and $A, B$, respectively. For this argument, we may even allow our choice of $A$ to lie outside of those usually allowed in the statistical model, as long as it lies in their (Zariski) closure. We therefore choose $A$ to be the identity, and $\boldsymbol{\pi}$ arbitrarily, so that $A'$ is also the identity, thus simplifying to considering

$$(8) \qquad B_1 = B_2 = B \otimes^{\text{row}} B \otimes^{\text{row}} \cdots \otimes^{\text{row}} B \qquad (k \text{ factors}).$$

It is now enough to show that $B_1$, as given in (8), has full row rank for some choice of stochastic $B$. We proceed very similarly to the proof of Lemma 13, but since a row tensor power occurs here rather than an arbitrary product, we must make some small changes to the argument.

Since nonzero rescalings of the rows of $B$ have no effect on the rank of $B_1$ in (8), we do not need to require that the row sums of $B$ are 1. So let $\mathbf{x} = (x_1, x_2, \ldots, x_\kappa)$ be a vector of distinct primes, and define $B$ by

$$
B = \begin{pmatrix}
1 & 1 & \cdots & 1 \\
x_1 & x_2 & \cdots & x_\kappa \\
x_1^2 & x_2^2 & \cdots & x_\kappa^2 \\
\vdots & \vdots & \ddots & \vdots \\
x_1^{r-1} & x_2^{r-1} & \cdots & x_\kappa^{r-1}
\end{pmatrix}.
$$

Let $\mathbf{y} = \mathbf{x} \otimes \mathbf{x} \otimes \cdots \otimes \mathbf{x}$ with $k$ factors. Then using notation from Lemma 13, $B_1$ will be the first $r$ rows of the Vandermonde matrix $W(\mathbf{y})$. To ensure $B_1$ has rank $r$, it is sufficient to ensure that $r$ of the entries of $\mathbf{y}$ are distinct, since then $B_1$ has a nonsingular Vandermonde submatrix of size $r$. The number of distinct entries of $\mathbf{y}$ is the number of distinct monomials of degree $k$ in the $x_i$, $1 \leq i \leq \kappa$. This number is $\binom{k+\kappa-1}{\kappa-1}$, so to ensure that generic $A, B$ lead to $B_1$ having full row rank, we ask that $k$ satisfy

$$\binom{k+\kappa-1}{\kappa-1} \geq r.$$

For fixed $\kappa$, the expression on the left of this inequality is an increasing and unbounded function of $k$, so this condition can be met for any $r, \kappa$.

Thus by Kruskal's theorem, from the joint distribution of $2k + 1$ consecutive variables of the HMM for generic $A, B$, we may determine $PB_1, PB_2, PB$, where $P$ is an unknown permutation.

Now to identify $A, B$ up to label swapping means to determine $\widetilde{A} = PAP^t$ and $\widetilde{B} = PB$ for some permutation $P$. As $\widetilde{B}$ has been found, we focus on $\widetilde{A}$. From (7) one finds

$$PB_2 = \widetilde{A}(\widetilde{B} \otimes^{\mathrm{row}} (\cdots \widetilde{A}(\widetilde{B} \otimes^{\mathrm{row}} (\widetilde{A}(\widetilde{B} \otimes^{\mathrm{row}} (\widetilde{A}\widetilde{B})))) \cdots)).$$

In this expression, $\widetilde{A}$ and $\widetilde{B}$ appear $k$ times. Since each row of $\widetilde{B}$ sums to 1, by appropriate summing of the columns of this matrix (marginalizing over the variable $X_{2k}$), we may determine a matrix $M$ given by a similar formula, but with only $k-1$ occurrences of $\widetilde{A}$ and $\widetilde{B}$. Then

$$PB_2 = \widetilde{A}(\widetilde{B} \otimes^{\mathrm{row}} M).$$

As $PB_2$ and $\widetilde{B} \otimes^{\mathrm{row}} M$ are known and generically of rank $r$, from this equation one can identify the matrix $\widetilde{A}$.

Thus the HMM parameters $A$ and $B$ are identifiable up to permutation of the states of the hidden variables. □

PROOF OF THEOREM 7. For the $n$ node model, with node set $V_n = \{v_k\}_{1 \leq k \leq n}$, denote (undirected) edges in the complete graph $K_n$ on $V_n$ by $(v_k, v_l) = (v_l, v_k)$ for $k \neq l$. We assume $0 < \pi_1 \leq \pi_2 < 1$, $p_{ij} \in [0,1]$ and $p_{11}, p_{12}, p_{22}$ are distinct.

Let $Z = (Z_1, Z_2, \ldots, Z_n)$ be the random variable, with state space $\{1, 2\}^n$, which describes the state of all nodes collectively. Ordering the elements of the state space in some way, we find the probabilities of the various states of $Z$ are given by the entries of a vector $\mathbf{v} \in \mathbb{R}^{2^n}$, all of which have the form $\pi_1^k \pi_2^{n-k}$. Observe for later use that no entries of $\mathbf{v}$ are zero, and the smallest and largest entries are $\pi_1^n$ and $\pi_2^n$, respectively (if $\pi_1 = \pi_2 = 1/2$, then all the entries of $\mathbf{v}$ are equal to $2^{-n}$).

Elements of the state space of $Z$ are specified by $\mathcal{I} = (i_1, i_2, \ldots, i_n) \in \{1, 2\}^n$, meaning $i_k$ is the state of $v_k$. An assignment of states to all edges in a graph $G \subseteq K_n$ will be represented by a subgraph $\mathcal{G} \subseteq G$ containing only those edges in state 1, in accord with the interpretation of edge states 0 and 1 as "absent" and "present." We refer to the probability of a particular state assignment to the edges of $G$ as the probability of *observing* the corresponding $\mathcal{G}$. We may think of any such $G$ as specifying a composite edge variable, whose states are represented by the subgraphs $\mathcal{G} \subseteq G$.

To relate the random graph model to the model of Kruskal's theorem, we must choose three observed variables and one hidden variable that reflect a conditional independence structure. The hidden variable will be $Z$ described above, indicating the state of some number of nodes $n$, to be chosen below. The observed variables will correspond to three pairwise edge-disjoint subgraphs $G_1, G_2, G_3$ of $K_n$. By choosing the $G_i$ to have no edges in common, we ensure that for $i \neq j$ observing any subgraph $\mathcal{G}_i$ of $G_i$ is independent of observing any subgraph $\mathcal{G}_j$ of $G_j$, conditioned on the state of $Z$. To meet the technical assumptions of Kruskal's theorem, we will also choose the $G_i$ so that the three matrices $B_i$ whose entries give probabilities of observing each subgraph $\mathcal{G}_i$ of $G_i$, conditioned on the state of $Z$, have full row rank. These matrices thus give conditional probabilities of observations marginalized over all edges not in $G_i$.

The construction of the $G_i$ proceeds in several steps. We begin by considering a small complete graph, and an associated matrix: for a set of 4 nodes, define a $2^4 \times 2^{\binom{4}{2}} = 16 \times 64$ matrix $A$, with rows indexed by assignments $\mathcal{I} \in \{1, 2\}^4$ of states to the nodes, columns indexed by all subgraphs $\mathcal{G}$ of $K_4$ and entries giving the probability of observing the subgraph conditioned on the state assignment of the nodes. Each entry of $A$ is thus a monomial in the $p_{ij}$ and $q_{ij} = 1 - p_{ij}$. Explicitly, if $\mathcal{I} = (i_1, i_2, i_3, i_4)$, and $e_{kl} \in \{0, 1\}$ is the state of edge $(v_k, v_l)$ in $\mathcal{G}$, the $(\mathcal{I}, \mathcal{G})$-entry of $A$ is

$$\prod_{1 \leq k < l \leq 4} p_{i_k i_l}^{e_{kl}} q_{i_k i_l}^{1 - e_{kl}}.$$

LEMMA 15. *For distinct $p_{11}, p_{12}, p_{22}$, the $16 \times 64$ matrix $A$ described above has full row rank.*

This lemma can be established by a rank computation with symbolic algebra software, so we omit a proof. One can also see, either through computation or reasoning, that the complete graph on fewer than 4 nodes fails to produce a matrix of full rank.

The next lemma shows we can find the 3 edge-disjoint subgraphs needed for the application of Kruskal's theorem. As the rest of the proof does not depend on the nodes having 2 states, we state the following lemma for an arbitrary number of node states. The more general formulation we prove here will be needed in a subsequent paper.

Let $r$ denote the number of node states and suppose we have found a number $m$ such that the $r^m \times 2^{\binom{m}{2}}$ matrix $A$ of probabilities of observations of subgraphs of the complete graph on $m$ nodes conditioned on node states has rank $r^m$. Lemma 15 establishes that for $r = 2$, we may take $m = 4$.

LEMMA 16. *Suppose for the $r$-node-state model, the number of nodes $m$ is such that the $r^m \times 2^{\binom{m}{2}}$ matrix $A$ of probabilities of observing subgraphs of $K_m$ conditioned on node state assignments has rank $r^m$. Then with $n = m^2$ there exist pairwise edge-disjoint subgraphs $G_1, G_2, G_3$ of $K_n$ such that for each $G_i$, the matrix $B_i$ of probabilities of observing subgraphs of $G_i$ conditioned on node state assignments has rank $r^n$.*

PROOF. We first describe the construction of the subgraphs $G_1, G_2, G_3$ of $K_n$. For each $G_i$, we partition the $m^2$ nodes into $m$ groups of size $m$ in a way to be described shortly. Then $G_i$ will be the union of the $m$ complete graphs on each partition set. Thus $G_i$ has $m\binom{m}{2}$ edges.

For conditional independence of observations of edges in $G_i$, from those in $G_j$ with $i \neq j$, we must ensure $G_i$ and $G_j$ have no edges in common. This requires only that a partition set of nodes leading to $G_i$ has at most one element in common with a partition set leading to $G_j$, if $i \neq j$. Labeling the nodes by $(i, j) \in \{1, \ldots, m\} \times \{1, \ldots, m\}$, we picture the nodes as lattice points in a square grid. We take as the partition leading to $G_1$ the rows of the grid, as the partition leading to $G_2$ the columns of the grid and as the partition leading to $G_3$ the diagonals. Explicitly, if $\mathcal{P}_i = \{V_j^i | j \in \{1, \ldots, m\}\}$ denotes the partition of the node set $V_{m^2}$ leading to $G_i$, then

$$V_j^1 = \{(j, i) | i \in \{1, \ldots, m\}\},$$
$$V_j^2 = \{(i, j) | i \in \{1, \ldots, m\}\},$$
$$V_j^3 = \{(i, i + j \bmod m) | i \in \{1, \ldots, m\}\}$$

and each $G_i$ is the union over $j \in \{1, \ldots, m\}$ of the complete graphs on node set $V_j^i$.

Now $B_i$, the matrix of conditional probabilities of observing all possible subgraphs of $G_i$, conditioned on node states, has $r^n$ rows indexed by composite states of all $n = m^2$ nodes, and $2^{m\binom{m}{2}}$ columns indexed by subgraphs of $G_i$. Observe that with an appropriate ordering of the rows and columns (which is dependent on $i$), $B_i$ has a block structure given by

$$(9) \qquad B_i = A \otimes A \otimes \cdots \otimes A \qquad (m \text{ factors}).$$

[Note that since $A$ is $r^m \times 2^{\binom{m}{2}}$, the tensor product on the right is $(r^m)^m \times (2^{\binom{m}{2}})^m$ which is $r^{m^2} \times 2^{m\binom{m}{2}}$, the size of $B_i$.] That $B_i$ is this tensor product is most easily seen by noting the partitioning of the $m^2$ nodes into $m$ disjoint sets $V_j^i$ gives rise to $m$ copies of the matrix $A$, one for each complete graph on a $V_j^i$. The row indices of $B_i$ are obtained by choosing an assignment of states to the nodes in $V_j^i$ for each $j$ independently, and the column indices by the union of independently-chosen subgraphs of the complete graphs on $V_j^i$ for each $j$. This independence in both rows and columns leads to the tensor decomposition of $B_i$.

Now since $A$ has full row rank, (9) implies that $B_i$ does as well. □

REMARK 1. For future work, we note that this lemma easily generalizes to graph models in which edges may be in any of $s$ states, with $s > 2$. In that case, the matrix $A$ is $r^m \times s^{\binom{m}{2}}$, and the columns of $A$ are no longer indexed by subgraphs of $K_m$, but rather by $s$-colorings of the edges of $K_m$.

To complete the proof of Theorem 7, we apply Corollary 2 for $\mathcal{M}(2^{m^2}; 2^{m\binom{m}{2}}, 2^{m\binom{m}{2}}, 2^{m\binom{m}{2}})$ to the parameter choice $\boldsymbol{\pi} = \mathbf{v}$, $M_i = B_i$, to find $\mathbf{v}$ and each $B_i$ is identifiable, up to row permutation. Thus here we do not apply the corollary to the full random graph model, but rather its marginalization over all edges not in $G_1 \cup G_2 \cup G_3$.

Suppose now that $\pi_1 \neq \pi_2$. Since $\pi_1^{m^2}, \pi_2^{m^2}$ are the smallest and largest entries of $\mathbf{v}$, respectively, we may determine $\pi_1, \pi_2$, as well as which of the rows of $B_i$ correspond to the having all nodes in state 1 or all in state 2. Summing appropriate entries of these rows, we obtain the probabilities $p_{11}$ and $p_{22}$ of observing a single edge conditioned on these node states. (This is simply a marginalization; sum the row entries corresponding to all subgraphs of $G_i$ that contain a fixed edge.) To find $p_{12}$, by consulting $\mathbf{v}$ we may choose one of the $n$ rows of $B_1$ which corresponds to node states with all nodes but one in state 1. By considering sums of row entries to obtain the conditional probability of observing a single edge, we can produce the numbers $p_{11}$, $p_{12}$. As $p_{11}$ is known, and $p_{12}$ is distinct from it, we thus determine $p_{12}$.

If $\pi_1 = \pi_2$, then we cannot immediately determine which rows of $B_i$ correspond to all nodes being in state 1 or in state 2. However, by marginalizing all rows to obtain the conditional probability of observing a single edge, we may determine the set of numbers $\{p_{11}, p_{12}, p_{22}\}$. With these in hand, we may then determine which two rows correspond to having all nodes in state 1 and all in state 2. This then uniquely determines which of the numbers is $p_{12}$, so everything is known up to label swapping. □

REMARK 2. While the above argument shows generic identifiability of the parameters of the 2-node state random graph model provided there are at least 16 nodes, a slightly more complicated argument, which we do not include here, can replace Lemma 16 to establish generic identifiability provided there are at least 10 nodes. Thus we make no claim to having determined the minimum number of nodes to ensure identifiability.

PROOF OF THEOREM 8.    We assume as usual that $Z$ is a latent random variable with distribution on $\{1, \ldots, r\}$ given by the vector $\boldsymbol{\pi}$, and $X = (X_1, \ldots, X_p)$ is the vector of observations such that conditional on $Z = i$, the variates $\{X_j\}_{1 \le j \le p}$ are independent, each $X_j$ having probability distribution $\mu_i^j$. We focus on 3 random variables at a time only, beginning first with $X_1, X_2, X_3$. The idea is to construct a binning of the random variables $X_1, X_2$ and $X_3$ using $\kappa_j - 1 \in \mathbb{N}$ cut points for $X_j$. For each $j = 1, 2, 3$, consider a partition of $\mathbb{R}$ into $\kappa_j$ consecutive intervals $\{I_j^k\}_{1 \le k \le \kappa_j}$, and consider the random variable $Y_j = (\mathbf{1}\{X_j \in I_j^1\}, \ldots, \mathbf{1}\{X_j \in I_j^{\kappa_j}\})$, where $\mathbf{1}\{A\}$ denotes the indicator function of set $A$. This is a finite random variable taking values in $\{0, 1\}^{\kappa_j}$ with at most one nonzero entry. We will show here that we can identify the proportions $\pi_i$ and the probability measures $\mu_i^j, 1 \le i \le r, 1 \le j \le 3$, relying only on the binned observed variables $\{Y_1, Y_2, Y_3\}$ for some well-chosen partitions of $\mathbb{R}$.

Consider for each $j = 1, 2, 3$, the matrices $M_j$ of size $r \times \kappa_j$ whose $i$th row is the distribution of $Y_j$ conditional on $Z = i$, namely the vector $[\mathbb{P}(X_j \in I_j^1 | Z = i), \ldots, \mathbb{P}(X_j \in I_j^{\kappa_j} | Z = i)]$. Introduce the matrix $\tilde{M}_1$ whose $i$th row is the $i$th row of $M_1$ multiplied by the value $\pi_i$. Note that the $M_j$'s are stochastic matrices. Moreover, the tensor product $[\tilde{M}_1, M_2, M_3]$ is the $\kappa_1 \times \kappa_2 \times \kappa_3$ table whose $(k_1, k_2, k_3)$ entry is the probability $\mathbb{P}((X_1, X_2, X_3) \in I_1^{k_1} \times I_2^{k_2} \times I_3^{k_3})$. This tensor is completely known as soon as the probability distribution (4) is given. Now we use Kruskal's result to prove that with knowledge of the tensor $[\tilde{M}_1, M_2, M_3]$, we can recover the parameters $\boldsymbol{\pi}$ and the stochastic matrices $M_j, j = 1, 2, 3$. If we can do so for general enough and well-chosen partitions $\{I_j^k\}_{1 \le k \le \kappa_j}$, then we will be able to recover the measures $\mu_i^j$, for $j = 1, 2, 3$ and $1 \le i \le r$.

We look for partitions $\{I^k\}_{1 \le k \le \kappa}$, with $\kappa \ge r$, such that the corresponding matrix $M$ has full row rank. Here we deliberately dropped the index $j = 1, 2, 3$. If we can construct these matrices with full row rank, then we get $I_1 + I_2 + I_3 = 3r \ge 2r + 2$ and Kruskal's result applies. As the partition $\{I^k\}_{1 \le k \le \kappa}$ is composed of consecutive intervals, the rows of the matrix $M$ are of the form $[F_i(u_1), F_i(u_2) - F_i(u_1), \ldots, F_i(u_{\kappa-1}) - F_i(u_{\kappa-2}), 1 - F_i(u_{\kappa-1})]$ for some real number cut points $u_1 < u_2 < \cdots < u_{\kappa-1}$. Replacing the $j$th column $C_j$ of $M$ by $C_j + C_{j-1}$ for consecutive $j$ from 2 to $\kappa$, we construct $M'$ with same rank as $M$ and whose $i$th row is $[F_i(u_1), F_i(u_2), \ldots, F_i(u_{\kappa-1}), 1]$. Now linear independence of the probability distributions $\{\mu_i\}_{1 \le i \le r}$ is equivalent to linear independence of the c.d.f.s $\{F_i\}_{1 \leq i \leq r}$. We need the following lemma.

LEMMA 17. *Let $\{F_i\}_{1 \leq i \leq r}$ be linearly independent functions on $\mathbb{R}$. Then there exists some $\kappa \in \mathbb{N}$ and real numbers $u_1 < u_2 < \cdots < u_{\kappa-1}$ such that the vectors*

$$\{(F_i(u_1), \ldots, F_i(u_{\kappa-1}), 1)\}_{1 \leq i \leq r}$$

*are linearly independent.*

PROOF. Let us consider a set of points $u_1 < u_2 < \cdots < u_m$ in $\mathbb{R}$ with $m \geq r$ and the matrix $A_m$ of size $r \times m$ whose $i$th row is $(F_i(u_1), \ldots, F_i(u_m), 1)$. Denote by $\mathcal{N}_m = \{\boldsymbol{\alpha} \in \mathbb{R}^r \mid \boldsymbol{\alpha} A_m = 0\}$, the left nullspace of the matrix $A_m$, and let $d_m$ be its dimension. If $d_m = 0$, then the matrix $A_m$ has full row rank and the proof is complete. Now if $d_m \geq 1$, choose a nonzero vector $\boldsymbol{\alpha} \in \mathcal{N}_m$. By linear independence of the $F_i$'s, we know that $\sum_{i=1}^{r} \alpha_i F_i$ is not the zero function, which means that there exists some $u_{m+1} \in \mathbb{R}$ such that $\sum_{i=1}^{r} \alpha_i F_i(u_{m+1}) \neq 0$. Up to a reordering of the $u$'s, we may assume $u_m < u_{m+1}$ and consider the matrix $A_{m+1}$ whose $i$th row is $[F_i(u_1), \ldots, F_i(u_{m+1}), 1]$ and whose left nullspace $\mathcal{N}_{m+1}$ has dimension $d_{m+1} < d_m$. Indeed, we have $\mathcal{N}_{m+1} \subset \mathcal{N}_m$ and by construction, the one-dimensional space spanned by the vector $\boldsymbol{\alpha}$ is not in $\mathcal{N}_{m+1}$. Repeating this construction a finite number of times, we find a matrix $A_\kappa$ with the desired properties. □

With this lemma, we have proved that the desired partition exists. Moreover, for any value $t \in \mathbb{R}$, we may, by increasing $\kappa$, include $t$ among the points $u_k$ without lowering the rank of the matrix. Thus we can construct partitions that involve any chosen cut point in such a fashion that Kruskal's result in the form of Corollary 2 will apply. That is, the vector $\boldsymbol{\pi}$ and the matrices $M_j$ may be recovered from the mixture $\mathbb{P}$, up to permutation of the rows. Moreover, summing up the first columns of the matrix $M_j$, up to the one corresponding to the chosen cut point $t$, we obtain the value of $F_i^j(t)$, $1 \leq i \leq r, j = 1, 2, 3$. To see that this enables one to recover the whole probability distribution $\mu_i^j$, up to label swapping on the $i$'s indexes, note that once we fix an ordering on the states of the hidden variable, the rows $(F_i(u_1), \ldots, F_i(u_\kappa))_{1 \leq i \leq r}$ are fixed and for each value of $t \in \mathbb{R}$, we associate to the $i$th row the value $F_i(t)$.

To conclude the proof, in the case of more than 3 variates, we repeat the same procedure with the random variables $X_1, X_2, X_4$. This enables us to recover the values of $\{\mu_i^1, \mu_i^2, \mu_i^4\}_{1 \leq i \leq r}$ up to a relabeling of the groups. As soon as the $\mu_i^1$ are linearly independent, they must be different, and using the two sets $\{\mu_i^1, \mu_i^2, \mu_i^3\}_{1 \leq i \leq r}$ and $\{\mu_i^1, \mu_i^2, \mu_i^4\}_{1 \leq i \leq r}$ which are each known only up to (different) label swappings, we can thus recover the set $\{\mu_i^1, \mu_i^2, \mu_i^3, \mu_i^4\}_{1 \le i \le r}$ up to a relabeling of the groups. Adding a new random variable at a time finally gives the result. □

PROOF OF THEOREM 9. In case of nonunidimensional blocks of independent components, we proceed much as in the proof of Theorem 8, but construct a binning into product intervals. For instance, if $X$ is two dimensional, we use $\kappa^2$ different bins, constructing $Y = (1\{X \in I^1 \times J^1\}, 1\{X \in I^1 \times J^2\}, \ldots, 1\{X \in I^\kappa \times J^\kappa\})$ where $\{J^k\}_{1 \le k \le \kappa}$ is a second partition of $\mathbb{R}$ into $\kappa \in \mathbb{N}$ consecutive intervals. This yields a matrix $M$ whose rows are of the form

$$(F_i(u_1, v_1), F_i(u_1, v_2) - F_i(u_1, v_1), \ldots, F_i(u_1, +\infty) - F_i(u_1, v_{\kappa-1}),$$
$$F_i(u_2, v_1) - F_i(u_1, v_1), F_i(u_2, v_2) - F_i(u_1, v_2)$$
$$- F_i(u_2, v_1) + F_i(u_1, v_1), \ldots,$$
$$F_i(u_2, +\infty) - F_i(u_1, +\infty) - F_i(u_2, v_{\kappa-1}) + F_i(u_1, v_{\kappa-1}), \ldots,$$
$$F_i(+\infty, v_1) - F_i(u_{\kappa-1}, v_1), F_i(+\infty, v_2) - F_i(u_{\kappa-1}, v_2)$$
$$- F_i(+\infty, v_1) + F_i(u_{\kappa-1}, v_1), \ldots,$$
$$1 - F_i(u_{\kappa-1}, +\infty) - F_i(+\infty, v_{\kappa-1}) + F_i(u_{\kappa-1}, v_{\kappa-1}))$$

for some real numbers $u_1 < u_2 < \cdots < u_{\kappa-1}$ and $v_1 < v_2 < \cdots < v_{\kappa-1}$. (To avoid cumbersome formulas, we only write the form of the matrix rows in the case $b = 2$.) This matrix has the same rank as $M'$ whose $i$th row is composed of the values $F_i(u_k, v_l)$ for $1 \le k, l \le \kappa$, using the convention $u_\kappa = v_\kappa = +\infty$. The equivalence between linear independence of the probability distributions and corresponding multidimensional c.d.f.'s remains valid.

Lemma 17 generalizes to the following.

LEMMA 18. *Let $\{F_i\}_{1 \le i \le r}$ be linearly independent functions on $\mathbb{R}^b$. There exists some $\kappa$, and $b$ collections of real numbers $u_1^i < u_2^i < \cdots < u_{\kappa-1}^i$, for $1 \le i \le b$, such that the $r$ row vectors composed of the values $\{F_i(u_{i_1}^1, \ldots, u_{i_b}^b) | i_1, \ldots, i_b \in \{1, \ldots, \kappa\}\}$, for $1 \le i \le r$ are linearly independent.*

The proof of this lemma is essentially the same as the proof of Lemma 17. The only difference with the previous setup is that now the construction of the desired set relies on addition of $b$ coordinates at a time, namely $t_1, \ldots, t_b \in \mathbb{R}$, which results in adding $\sum_{j=0}^{b-1} \binom{b}{j} \kappa^j$ columns in the matrix.

To complete the argument establishing Theorem 9, we may again include any point $(t_1, \ldots, t_b) \in \mathbb{R}^b$ among the $u_k$'s without changing the row ranks of the matrices to which we apply Kruskal's theorem. Thus we may recover the values $F_1(t_1, \ldots, t_b), \ldots, F_r(t_1, \ldots, t_b)$, and conclude the proof in the same way as the last theorem.    □

PROOF OF LEMMA 10.    Suppose the probability distribution

$$\mathbb{P}(X_1, X_2) = \sum_{i=1}^{r} \pi_i \mu_i^1(X_1) \mu_i^2(X_2)$$

has rank $r$. Then for $k = 1, 2$, the sets $\{\mu_i^k\}_{1 \leq i \leq r}$ must be independent, since any dependency relation would allow $\mathbb{P}$ to be expressed as a sum of fewer products.

Conversely, suppose for $k = 1, 2$, the measures $\{\mu_i^k\}_{1 \leq i \leq r}$ are independent. The corresponding sets of c.d.f.s $\{F_i^k\}_{1 \leq i \leq r}$ are also independent, and thus we may choose collections of points $\{t_j^k\}_{1 \leq j \leq r}$ such that the $r \times r$ matrices $M_k$ whose $i, j$-entries are $F_i^k(t_j^k)$ have full rank. Then with $F$ denoting the c.d.f. for $\mathbb{P}$, the matrix $N$ with entries $F(t_i^1, t_j^2)$ can be expressed as

$$N = M_1^T \operatorname{diag}(\boldsymbol{\pi}) M_2$$

and therefore has full rank. But if the rank of $\mathbb{P}$ were less than $r$, a similar factorization arising from the expression of $\mathbb{P}$ using fewer than $r$ summands shows that $N$ has rank smaller than $r$. Thus the rank of $\mathbb{P}$ is at least $r$, and since the given form of $\mathbb{P}$ shows the rank is at most $r$, it has rank exactly $r$. □

**Acknowledgments.**    The authors thank the Isaac Newton Institute for the Mathematical Sciences and the Statistical and Applied Mathematical Sciences Institute for their support during residencies in which some of this work was undertaken.

E. S. ALLMAN
J. A. RHODES
DEPARTMENT OF MATHEMATICS AND STATISTICS
UNIVERSITY OF ALASKA, FAIRBANKS
FAIRBANKS, ALASKA 99775
USA
E-MAIL: e.allman@uaf.edu
            j.rhodes@uaf.edu

C. MATIAS
CNRS UMR 8071
LABORATOIRE STATISTIQUE ET GÉNOME
523, PLACE DES TERRASSES DE L'AGORA
91 000 ÉVRY
FRANCE
E-MAIL: catherine.matias@genopole.cnrs.fr